\input amstex.tex
\documentstyle{amsppt}
\input xy
\xyoption{all}

\magnification=1200
\hsize=150truemm
\vsize=224.4truemm
\hoffset=4.8truemm
\voffset=12truemm

\NoBlackBoxes
\NoRunningHeads

\def\Square{\rlap{$\sqcup$}$\sqcap$}
\def\cqfd {\quad \hglue 7pt\par\vskip-\baselineskip\vskip-\parskip
{\rightline{\Square}}}

\define\T{{\Cal T}}
\define\E{{\Cal E}}
\redefine\H{{\Cal H}}
\define\s{{\Cal S}}
\define\cc{{\Cal C}}
\define\R{{\bold R}}
\define\Z{{\bold Z}}
\define\mi{^{-1}}
\let\thm\proclaim
\let\fthm\endproclaim
\let\inc\subset 
\let\ds\displaystyle
\let\bo\partial
\define\Aut{\text{\rm Aut}}
\define\Out{\text{\rm Out}}
\define\mcg{MCG}
\define\rmcg{MCG^{\bo}}
\define\autr{\Aut_\H^{\bo}}
\define\outr{\rmcg}
\define\auta{\Aut_\H}
\define\outa{\mcg}

\newcount\secno
\newcount\subsecno
\newcount\stno
\global\subsecno=1

\define\sta{\the\secno.\the\stno
\global\advance\stno by 1}

\define\sect{\global\advance\secno by 1\global\subsecno=1\global\stno=1\
\the\secno. }

\def\nom#1{\edef#1{\the\secno.\the\stno}}

\newcount\refno
\global\refno=0
\def\nextref#1{\global\advance\refno by 1\xdef#1{\the\refno}}
\def\bref {\ref\global\advance\refno by 1\key{\the\refno}}


\nextref\BJ
\nextref\BF
\nextref\BH
\nextref\Bow
\nextref\CT
\nextref\Fuj
\nextref\HM
\nextref\IMc
\nextref\LL
\nextref\MNS
\nextref\Pau
\nextref\Pet
\nextref\Ripsel
\nextref\Sela
\nextref\Shor

\topmatter

\abstract Using the canonical JSJ splitting, we describe the outer automorphism
group $\Out(G)$ of a one-ended word hyperbolic group $G$. In particular, we
discuss to what extent $\Out(G)$ is virtually a direct product of mapping class
groups   and a free abelian group, and  we determine for which groups
$\Out(G)$ is infinite. We also show that there are only finitely many conjugacy
classes of torsion elements in $\Out(G)$, for $G$ any torsion-free hyperbolic
group.

More generally, let $\Gamma $ be a finite graph of groups decomposition of an
arbitrary group $G$ such that edge groups $G_e$ are rigid (i.e\. $\Out(G_e)$ is
finite). We describe   the group of automorphisms of $G$
preserving   $\Gamma $, by comparing it to direct products of
suitably defined mapping class groups of vertex groups.

\endabstract

\title Automorphisms of   hyperbolic groups and graphs of groups
 \endtitle

\author  Gilbert Levitt
\endauthor

\endtopmatter

\document

\head  {\sect Introduction and statement of results  }\endhead

The basic idea for understanding the outer automorphism group $\Out(G)$ of a
torsion-free one-ended word
hyperbolic group $G$ is due to Rips and Sela [\Ripsel, \Sela]. Since $G$ has a
canonical JSJ splitting, automorphisms of $G$ may be seen on this splitting:
$\Out(G)$ is virtually generated by Dehn twists along edges of the splitting,
together with mapping class groups of punctured surfaces (automorphism groups of
``quadratically hanging''  subgroups).

Giving a precise description of $\Out(G)$, however, requires a
careful study of these generators. This is what we do here, both for hyperbolic
groups and for splittings of arbitrary groups.

Let us  first illustrate this on a simple example. Let $\Sigma  $ be a $2$-torus
with one open disc removed.
Let $G$ be
the fundamental group of the
$2$-complex $K $ obtained by glueing three copies   $\Sigma _1$, $\Sigma
_2$, $\Sigma _3$ of $\Sigma $ along their boundary circle $C$. All
automorphisms of $G$ are realized by homeomorphisms of $K $. Since
homeomorphisms may permute the $\Sigma _i$'s, or reverse orientation, we
restrict
to homeomorphisms mapping each $\Sigma _i$ to itself in an
orientation-preserving
way; this amounts to passing to a finite index subgroup $\Out_0(G)\inc\Out(G)$.

Restricting homeomorphisms to each $\Sigma _i$ gives an epimorphism $\ds\rho
:\Out_0(G)\to\prod_{i=1}^3\mcg(\Sigma _i),$
where  $\mcg(\Sigma _i)\simeq \Out(F_2)\simeq SL(2,\Z)$ is the mapping
class group
of the punctured torus, i.e\. the group of isotopy classes of
orientation-preserving homeomorphisms of
$\Sigma _i$. The kernel of $\rho $ is generated by $D_1$, $D_2$, $D_3$, with
$D_i$ the (automorphism induced by a)  Dehn twist in an annulus parallel to the
boundary in $\Sigma _i$. The product of the $D_i$'s being isotopic to the
identity, $\ker\rho $ is isomorphic to $\Z^2$.

The epimorphism $\rho $ does not have a section (basically because the map
$\Aut(F_2)\to \Out(F_2)$ has none).
In
geometric terms, there is no canonical way of glueing elements of
$\mcg(\Sigma _i)$ along
$C$, because we consider homeomorphisms of
$\Sigma _i$ up to an isotopy which is free on the boundary. We therefore
introduce
$\rmcg(\Sigma _i)$, the group of
homeomorphisms of
$\Sigma _i$ fixing the boundary pointwise, up to isotopy
{relative to the boundary}. This group is a nontrivial central extension of
$\mcg(\Sigma _i)$ by $\Z$ (generated by the twist $D_i$).

Now we can  glue   along $C$, thus getting an
epimorphism
$ \lambda :\prod_{i=1}^3\rmcg(\Sigma _i)\to\Out_0(G)$,
whose kernel is isomorphic to $\Z$ (corresponding to the relation
$D_1D_2D_3=1$).

The discussion of this example is summed up by the following commutative
triangle
of nontrivial central extensions:

$$
\xymatrix
{\Z\ \ar@{^{(}->} [dr]
\\ \Z^3\ \ar @{^{(}->}[r]&{\prod\rmcg(\Sigma _i)}\ar @{>>}[rr]^{\pi
}\ar@{>>} [rd]
_{\lambda }
&&{\prod\mcg(\Sigma _i)}\\ &&{\Out_0(G)}\ar@{>>} [ur]_{\rho }\\
&\Z^2\ \ar@{^{(}->} [ur]
}
$$

Now consider an arbitrary one-ended hyperbolic group $G$ (assumed to be
torsion-free for the moment). We use the JSJ splitting of $G$
constructed by Bowditch [\Bow], because it is completely canonical and thus
invariant  under all automorphisms of $G$.

It is a finite graph of groups $\Gamma $, with $\pi
_1(\Gamma )=G$. There are three types of vertices. Let $V_1$ denote the set of
vertices
of
$\Gamma
$ with cyclic vertex group (we use the same 1, 2, 3
numbering   as in [\Bow]). Let
$V_2$  be the set of vertices $v$
 corresponding to ``MHF subgroups'': the vertex group
$G_v$ is the fundamental group of a compact surface $\Sigma _v$, with edge
groups
being fundamental groups of boundary components (in the example studied above
$V_1$ has one element, corresponding to $C$, and $V_2$ has three
elements, corresponding to the $\Sigma _i$'s).  Let
$V_3$ consist of the  remaining vertices of $\Gamma $, which we call rigid. Let
 $\E$ be the set of (unoriented) edges of
$\Gamma $.  Every edge   joins a vertex of $V_1$ to
a vertex of $V_2\cup V_3$, leading to $\E=\E_2\cup \E_3$.

Given a compact surface $\Sigma $, we denote by $\mcg(\Sigma )$ the (pure)
mapping
class group of the  surface with punctures, that is the group of
homeomorphisms of
$\Sigma
$ fixing the boundary pointwise, modulo those that are freely isotopic to
the identity (the boundary may ``turn'' during the isotopy).
We also define $\rmcg(\Sigma )$, the mapping class group of the surface with
boundary,   dividing only by homeomorphisms isotopic to the identity
relatively to the boundary. It is a central extension of $\mcg(\Sigma )$ by
$\Z^b$,
where
$b$ is the number of boundary components of $\Sigma $.

As suggested in [\Bow], the $\Aut(G)$-invariant JSJ splitting allows a simple
proof of the following result due to Sela.

\nom\prem
\thm{Theorem \sta {} ([\Sela, Theorem 1.9])} Let $G$ be a torsion-free one-ended
hyperbolic group. There is an exact sequence
$$1\to\Z^n\to\Out'(G)\to\prod_{v\in V_2}\mcg(\Sigma _v)\to1,
$$
where $\Out'(G)$ has finite index in $\Out(G)$, and $n=|\E|-|V_1|$ is
``the number of edges in the essential
JSJ decomposition of $G$''.
\fthm

This central extension, however, is not trivial in general, so one cannot claim,
as in [\Sela], that $\Out(G)$ is virtually a direct product of mapping class
groups and a free abelian group. The next result  leads to such a product
structure, but only under an additional hypothesis.

 \nom\deux

\thm{Theorem \sta {}  } Let $G$ be a torsion-free one-ended hyperbolic
group.  Let $r$ be the
number of vertices of
$V_1$ connected only to vertices of $V_2$, and $q=|\E_3|-|V_1|+r$. The group
$\Out(G)$ is virtually a direct product $\Z^q\times M$, where $M$ is the
quotient
of $\prod_{v\in V_2}\rmcg(\Sigma _v)$ by a  central  subgroup
isomorphic to
$\Z^r$.
\fthm

\thm{Corollary \sta} If every cyclic vertex of $\Gamma $ is connected to at
least
one rigid vertex, then $\Out(G)$ is virtually a direct product of a free abelian
group and mapping class groups $\rmcg(\Sigma _v)$ of surfaces with boundary.
\fthm

There are statements similar to Theorems \prem{} and \deux{} if $G$ is
allowed to
have torsion, they will be given in section 5. The main differences are that
$\mcg(\Sigma _v)$ and $\rmcg(\Sigma _v)$ are now mapping class groups of
$2$-orbifolds, and only edges of $\Gamma $ whose group has infinite center
must be
considered. Our study allows us to prove the following result (see section 6).
It was conjectured in [\MNS] (we thank G\. Swarup for suggesting this
application).

\nom\conj
\thm{Theorem \sta} Let $G$ be a one-ended hyperbolic group. Then $\Out(G)$ is
infinite if and only if $G$ splits  over a virtually cyclic
subgroup with infinite center, as an arbitrary HNN extension or as an amalgam of
groups with finite center.
\fthm

Paulin's theorem [\Pau], together with Rips theory [\BF], provides a
splitting over
a  virtually cyclic
subgroup, but with no control on the center.
Torsion-free hyperbolic groups with
infinitely many ends always have $\Out(G)$ infinite. Allowing torsion makes
things
much more complicated, see [\Pet] for the case of virtually free groups.

We also prove (see section 7):

\nom\nombfi
\thm{Theorem \sta} If $G$ is a  torsion-free  hyperbolic group, there are only
finitely many conjugacy classes of torsion elements in
$\Out(G)$.
\fthm

This will be used in [\LL] to give a proof of Shor's theorem on fixed subgroups
[\Shor].

\medskip

Now let $\Gamma $ be a finite graph of groups decomposition of an arbitrary
group
$G$.  We consider the group $\Out^\Gamma (G)\inc\Out(G)$ consisting of
automorphisms preserving $\Gamma $. This group was thoroughly studied   by
Bass-Jiang
[\BJ]. Here we focus on  a subgroup   $\Out_1^\Gamma (G)\inc \Out^\Gamma (G)$,
which we squeeze between products of automorphism groups associated to vertex
groups.

Let $V$ be the vertex set of $\Gamma $. For $v\in V$ we consider the vertex
group
$G_v$, with incident edge groups $G_e$. We let $\mcg(G_v)\inc\Out(G_v)$ be the
group of automorphisms which act on each $G_e$ as conjugation by some $g_e\in
G_v$. We also define an extension $\rmcg(G_v)$, by keeping track of the elements
$g_e$ (see section 4). These groups are an algebraic generalization of the
groups
$\mcg(\Sigma _v)$ and
$\rmcg(\Sigma _v)$ considered above.

\thm {Theorem \sta} There are natural epimorphisms $\lambda :
\prod_{v\in V}\rmcg(G _v)\to\Out^\Gamma _1(G)$ and
$\rho _1:\Out^\Gamma _1(G)\to\prod_{v\in V}\mcg(G _v)$, with $\Out^\Gamma
_1(G)\inc \Out^\Gamma (G)$ and $\rho _1\circ \lambda $ equal to the natural
projection.

If
$\Out(G_e)$ is finite for every edge group
$G_e$, then  $\Out^\Gamma _1(G)$ has finite index in
 $\Out^\Gamma (G)$.
\fthm

The kernels of $\lambda $ and $\rho _1$ may be described
explicitly. In particular, our main technical result (Proposition 3.1) is a
presentation for
$\T=\ker\rho _1$, the group of twists around edges of $\Gamma $.

\head  {\sect Automorphisms of graphs of groups  }\endhead

Most   results in this section follow from  [\BJ], but our exposition will
be  self-contained.

\subhead  {  Generalities  }\endsubhead

The following notation will be used throughout the paper. If $H$ is a
subgroup of
a group $G$, then $Z_G(H)$ and $N_G(H)$ denote the {\it centralizer\/} and
{\it normalizer\/} of $H$ in $G$, and $Z(H)$ is the {\it center\/} of $H$. For
$m\in G$, denote by
$i_m$ the inner automorphism $g\mapsto mgm\mi$. If $\alpha \in\Aut(G)$, let
$\hat\alpha $ be its class in $\Out(G)$. We say that $\alpha $ {\it
represents\/}
$\hat\alpha $.

Let $G$ be the fundamental group of a finite  graph of groups $\Gamma
$, with  vertex set $V$, vertex groups $G_v$, and edge groups $G_e$.  Let
$E$ be the set of {\it oriented edges\/} $e$ of $\Gamma $, and $E_v$   the set
of    edges $e$ with origin $o(e)=v$. Let $\E$ be
the set of unoriented edges. The incident edge groups
$G_e$, for $e\in E_v$, and their conjugates in $G_v$, will be called   {\it
peripheral subgroups\/} of
$G_v$.

Let $T$
be the  Bass-Serre $G$@-tree associated to $\Gamma $. We simply  write $g$
for the
automorphism of $T$ associated to $g\in G$, and $G_x$ for the stabilizer of
$x\in
T$. Choosing a fundamental domain, we identify each
 vertex
$v$ of $
\Gamma
$ with a vertex of $T$ whose stabilizer is $G_v $. Similarly, we often write $e$
for an edge of $T$ with stabilizer $G_e$.

We assume  that $\Gamma $ is {\it minimal\/}:   $G_e $ is a proper
subgroup of $G_v$ for every terminal vertex $v$ (equivalently, $\pi
_1(\Gamma ')$
is a proper subgroup of $\pi _1(\Gamma )$ for every proper connected subgraph
$\Gamma '$). We also assume that $\Gamma $ is not a {\it mapping torus\/}:
we say
that $\Gamma
$ is a  mapping torus if
$\Gamma
$ is topologically a circle and all inclusions from edge group to vertex
group are
onto (equivalently, $T=\R$ and $G$ acts by translations); if $\Gamma $ is a
mapping torus, then
$G$ is   a semi-direct product $G_v\rtimes_\varphi \Z$, with $\varphi
\in\Aut(G_v)$.

These assumptions imply that the centralizer of the image of $G$ in $\Aut(T)$ is
trivial. In particular, $Z(G)$ acts trivially on $T$.

 Let $\Aut^\Gamma (G)\inc\Aut(G)$ and $\Out^\Gamma (G)\inc\Out(G)$ be the
groups of
{\it automorphisms preserving $\Gamma $.\/} Topologically, $G$ is the
fundamental
group of a graph of spaces, and automorphisms in  $\Out^\Gamma (G)$ are
induced by
homeomorphisms preserving such a structure. Algebraically, the condition is best
expressed in terms of
$T$: an  automorphism
$\alpha
\in\Aut(G)$ belongs to
$\Aut^\Gamma (G)$ if and only if there exists an automorphism $H_\alpha $ of $T
$ conjugating the action of $G$ with the action twisted by $\alpha $, in
the sense
that
$\alpha (g)H_\alpha =H_\alpha  g$ for every $g\in G$ (if $G$ fixes no end
of $T$,
this is the same as saying that the length function satisfies
$\ell\circ\alpha =\ell$).

With our assumptions, $H_\alpha $ is unique and $\alpha \mapsto H_\alpha $
defines
an action of $\Aut^\Gamma (G) $ on $T$, with $i_m$ acting as $m$ for $m\in
G$. In
the situation of [\Bow],
this action is induced by the action of $\Aut(G)$ on $\bo G$.

\subhead  {  The homomorphism $\rho :\Out^\Gamma _0(G)\to\prod_{v\in
V}\Out(G_v)$
}\endsubhead

The map $H_\alpha $ induces an automorphism of the finite graph $\Gamma
=T/G$, and
we let
$\Out^\Gamma _0(G)\inc\Out^\Gamma (G)$ be the finite index subgroup
consisting of
automorphisms acting as the identity on   $\Gamma $.

Note that for $x\in T$ the stabilizer of $H_\alpha x$ is $\alpha (G_x)$, and
points in the same $G$-orbit have conjugate stabilizers.
If
$\hat \alpha \in\Out^\Gamma _0(G)$ and $v\in V$, then $\alpha (G_v)$ is a
conjugate
of
$G_v$ and we can define $\rho _v(\hat\alpha )\in\Out(G_v)$. More precisely,
choose $m\in G$ such
that
$H_\alpha v=mv$ (recall that we identify $v$ with a vertex of $T$). Then
$\alpha (G_v)=G_{mv}=mG_vm\mi$, so that
 $\beta =i_{ m\mi}\alpha  $ maps $G_v$ to itself. It is  easily checked that the
class of
$\beta
$ in
$\Out(G_v)$   depends only   on $\hat\alpha $, and we obtain a homomorphism
$\rho _v:\Out^\Gamma _0(G)\to\Out(G_v)$.

We also define $\rho =\prod
\rho _v:\Out^\Gamma _0(G)\to\prod_{v\in V}\Out(G_v)$. In the rest of this
section
we study the image and kernel of $\rho $.

\subhead  { Extending vertex automorphisms  }\endsubhead

Elements in the image of $\rho _v$ preserve the peripheral structure of $G_v$
(if $e$ is an edge of $T$ containing $v$, then
$H_\beta e=g_ee$ for some $g_e\in G_v$; the stabilizer of $H_\beta e$ is
$\beta(G_e)$, and also $g_eG_eg_e\mi$). The converse is not necessarily true: if
$G=\Z*_{2\Z}H$, the nontrivial automorphism of $\Z$ does not always extends
to $G$.

Define the mapping class group $\mcg(G_v)\inc \Out(G_v)$ by restricting to
automorphisms which act on each edge group $G_e$, $e\in E_v$, as conjugation by
some
$g_e\in G_v$.

\nom\debut
 \thm{Proposition \sta}
The product
$\prod_{v\in V}\mcg(G_v)$ is contained in the image of  $\rho $.
\fthm

\demo{Proof}
Given $\beta \in\Aut(G_v)$ with $\hat\beta \in\mcg(G_v)$, we   extend
it to
$\alpha
\in\Aut(G)$ with
$\rho _v(\hat\alpha )=\hat\beta $ and $\rho
_w(\hat\alpha ) =1$ for $w\neq v$.
For each
$e\in E_v$, fix
$g_e\in G_v$ such that
$\beta (g)=g_egg_e\mi$ for $g\in G_e$.

For concreteness we start with the
elementary cases.
   If
$G=G_v*_{G_e}H$, we define
$\alpha $ on $H$ as conjugation by $g_e$. If $G$ is an HNN-extension $<G_v,t\mid
t\varphi (g)t\mi=g\ \text{if}\ g\in G_e>$, we have $\beta (g)=g_egg_e\mi$ and
$\beta (\varphi (g))=g_f\varphi (g)g_f\mi$ for $g\in G_e$, and  we define
$\alpha
(t)=g_et g_f\mi$.

In the general case, recall that $G$ is a subgroup of a path group $\pi
(\Gamma )$   generated by the $G_v$'s and the edges of $\Gamma $ (see e.g\.
[\BJ]). We define
$\alpha
$ on the generators of
$\pi (\Gamma )$ in the following way: it equals $\beta $ on
$G_v$, it is the identity on
the other vertex groups and on edges not incident to $v$, and an edge $e$ with
origin $v$ is mapped to $g_ee$ (to $g_eeg_{\bar e}\mi$ if the opposite edge
$\bar
e$ also has origin $v$).  One checks that this is compatible with the relations
defining
$\pi (\Gamma )$, and   that the subgroup
$G$
  is invariant.
\cqfd\enddemo

This will be referred to as the {\it
extension construction\/}.
It  is simply  an
   algebraic translation of the following topological observation: if $K_0$ is a
subcomplex of codimension $0$  of a complex
$K $, every self-homeomorphism of $K _0$ which equals the identity on
$\bo K_0$ has a canonical extension to $K $.

\subhead  {  Twists and bitwists
}\endsubhead

 The extension construction depends on the choice  of  the elements $g_e$.
Each
$g_e$  may be  right-multiplied by an element $z$ of the centralizer
$Z_{G_v}(G_e)$. This changes $\alpha $ by   $D_z$,  the
extension of
$\beta =1$ relative to the choice $g_e=z$ (and $g_f=1$ for $f\neq e$).

We call $D_z$ the {\it twist by $z$ around $e$ (near $v$)\/}.  It may also be
defined, more directly, as  follows. Let $e\in E_v$ and $z\in
Z_{G_v}(G_e)$. If $e$ separates
$\Gamma
$, we have an amalgam and
$D_z$ acts as the identity on the factor containing $G_v$ and as
conjugation by $z$
on the other. In the case of an HNN-extension, $D_z$ maps the stable letter
$t$ to
$zt$. If $z$ is a generator of $G_e\simeq\Z$, then $D_z$ is induced by a
Dehn twist
on an annulus.

We denote by $\T$ the
subgroup of $\Out^\Gamma _0(G)$ generated by all twists $\hat D_z$.
Extensions of a
given $\hat\beta $ associated to different choices of the $g_e$'s differ by an
element of $\T$.

We will see that $\T$ has finite index in $\ker\rho $ when all groups
$\Out(G_e)$
are finite. In general, though, describing $\ker\rho $ requires
   bitwists, which we now define.

Let $e$ be an edge of $\Gamma $ with  endpoints $v,w$ (possibly equal). Suppose
that
$z\in G_v$ normalizes the image of $G_e$, that $z'\in G_w$ normalizes the
image of
$G_e$, and that
$z,z'$ have the same action on $G_e$ (if
$v\neq w$, this means that $z\mi z'$ centralizes $G_e$). We define the {\it
bitwist\/}
$D_{z,z'}$ around $e$ as conjugation by $z$ on one factor and by
$z'$ on the other in the case of an amalgam, as mapping $t$ to $z  tz'{}\mi$
in the HNN case. For instance, amalgamate two Klein bottle groups $<t_i a
t_i\mi=a\mi>$ along the cyclic group $<a>$. Then $D_{t_1,t_2}$ fixes
$t_1,t_2$ and
maps $a$ to
$a\mi$.

Of course twists are bitwists, and (bi)twists
around distinct edges commute in
$\Out(G)$. Also note that
bitwists belong to $\T$ when edge groups are trivial, and   when $G$ is a
torsion-free hyperbolic group and edge groups are cyclic (because no nontrivial
element of $G$ is conjugate to its inverse).

\subhead The kernel of $\rho $ \endsubhead

\nom\noyau
\thm{Proposition \sta} The kernel of $\rho
$ is generated by bitwists around the edges of $\Gamma $.
\fthm

\demo{Proof} Fix a maximal subtree
$\Gamma _0\inc\Gamma
$, which we identify with a lift to $T$.
Let $e =vw$
be an edge of $\Gamma _0$. Any $\hat\alpha \in\ker\rho $ has a representative
$\alpha $,   equal to the identity on $G_v$, with $H_\alpha v=v$. The image
of the
edge
$e
$ by $H_\alpha $  is an edge in the same $G$-orbit, which we can write
$z\mi e$ for some
$z\in G_v$. Note that $z$ normalizes $G_e$ since the stabilizer of $H_\alpha
e=z\mi e$ is both $\alpha (G_e)=G_e$ and $z\mi G_ez$. Now $\alpha (G_w)$ is the
stabilizer of $H_\alpha w=z\mi w$, equal to
$z\mi G_wz$. Since $\rho _w(\hat\alpha )$ is trivial,
there exists $z'\in G_w$ such that $\alpha (g)= z\mi z' gz'{}\mi z$ for $g\in
G_w$.
We conclude that $\hat\alpha \circ \hat D_{z,z'}\mi$ has a representative
$\alpha
'$ inducing the identity on $<G_v,G_w>$, with $H_{\alpha '}$ equal to the
identity
on $e$.

Repeating this argument, we can compose any $\hat\alpha \in \ker\rho $ by
bitwists and obtain $\hat\alpha '$ with a representative $\alpha '$ equal to the
identity on the fundamental group of $\Gamma _0$, with $H_{\alpha '}$ equal
to the
identity on $\Gamma _0$. If $\Gamma _0=\Gamma $ we are done since
$\hat\alpha $ is
a product of bitwists. Otherwise we apply a similar argument to the
HNN generators.
\cqfd\enddemo

\subhead  {  Rigid edge groups
}\endsubhead

\nom\indice
\thm{Proposition \sta} Suppose that $\Out(G_e)$ is finite for every edge group
$G_e$ of $\Gamma $. Then
$\prod_{v\in V}\mcg(G_v)$ has finite index in the image of   $\rho
$, and the group of twists $\T$ has
finite index in $\ker\rho $.
\fthm

\demo{Proof} For the first assertion, it suffices to show that $\mcg(G_v)$ has
finite index in the image of $\rho _v$ for every $v\in V$. Since   any
$\varphi $
in the image of $\rho _v$ preserves the peripheral structure of $G _v$, one can
associate to $\varphi $ and
$e\in E_v$ a right   coset of  $\Out(G_e)$ modulo the image of
$N_{G_v}(G_e)$. If $\varphi $ and $\varphi '$ define the same coset for
every $e\in
E_v$, then $\varphi $ and $\varphi '$ belong to the same coset modulo
$\mcg(G_v)$.

The second assertion follows from the fact that $\hat D_{z,z'}$ belongs to
$\T$ if
  $z$ and $z'$ act on $G_e$  by an inner automorphism.
\cqfd\enddemo

\medskip

\head  {\sect The group of twists $\T$ }\endhead

We view $\T$ as a quotient
of
$\ds\prod_{v\in
V}\biggl({\prod_{e\in E_v} Z_{G_v}(G_e)}\biggr)=\prod_{e\in E}
Z_{G_{o(e)}}(G_e)$.
We denote by $\theta $ the quotient map. Generators for $\ker \theta $ will be
given by {\it  vertex and edge relations\/}, coming from centers of vertex
groups and edge groups: twisting by
$z\in Z(G_v)$ around all edges with origin $v$, or by $z\in Z(G_e)$ at both ends
of $e$, produces an inner automorphism.

More formally, we embed
$ Z(G_v)$  into $\prod_{e\in E} Z_{G_{o(e)}}(G_e)$ by embedding it
diagonally into
${\prod_{e\in E_v} Z_{G_v}(G_e)}$. Given an edge $e$ with endpoints $v,w$,
we have
embeddings of
$G_e$ into $G_v$ and $G_w$, therefore a diagonal embedding of   $ Z(G_e)$
 into
$Z_{G_v}(G_e)\times Z_{G_w}(G_e)$. Putting all these embeddings together gives a
map $$j:\prod_{v\in V}Z(G_v)\times
 {\prod_{e\in \E }
Z (G_e)}\to\prod_{e\in E} Z_{G_{o(e)}}(G_e)$$ (recall that $\E$ is
the set  of unoriented edges  and $E$ is the set    of oriented edges).
Note that
only two terms contribute to a given factor
$Z_{G_v}(G_e)$, namely $Z(G_v)$ and $Z(G_e)$. The image of
$j$ is central and contained in $\ker \theta $.

\nom\relat
\thm{Proposition \sta} Assume that $\Gamma $ is minimal and is not a
mapping torus.
The group of twists
$\T$ is the quotient of
$\ds\prod_{e\in E} Z_{G_{o(e)}}(G_e)$ by the vertex and edge relations: there is
an exact sequence
$$\{1\}\to Z(G)\overset i \to\to
\prod_{v\in V}Z(G_v)\times
 {\prod_{e\in \E }
Z (G_e)}
\overset j \to\to
\prod_{e\in E} Z_{G_{o(e)}}(G_e)
\overset \theta \to\to
\T\to\{1\}.
$$

\fthm

\demo{Proof}  Elements of $Z(G)$ act trivially on $T$, hence belong to every
$Z(G_v)$ and $Z(G_e)$. We define $i$ by mapping $n\in Z(G)$ to $n$ in  every
$Z(G_v)$ and $n\mi$ in every $Z(G_e)$. It is easy to check that the image of
$i$ equals $\ker j$.

Given
$\underline z=(z_e)_{e\in E}$ with
$\theta (\underline z)$ trivial,  we must now show that
$\underline z$ belongs to the image $J$ of $j$.  The proof is by induction
on the
number of edges of $\Gamma $, we divide it into four steps.

{\it $\bullet$  First step: graphs
with one edge.}

First suppose that $\Gamma $ is an amalgam $G=G_1*_{G_e}G_2$. Consider two
elements $z_i\in Z_{G_i}(G_e)$ such that the
automorphism equal to conjugation by $z_i$ on $G_i$ is inner. There exists $m\in
G$ such that $u_i=m\mi z_i$ commutes with $G_i$. Since $\Gamma $ is minimal we
have $N_G( G_i)=G_i$, and therefore
$u_i\in Z(G_i)$. Furthermore $m=z_1u_1\mi=z_2u_2\mi\in Z(G_e)$, showing that
$(z_1,z_2)$ belongs to   $J$ since
$z_i=mu_i$ with $m\in Z(G_e)$ and $u_i\in Z(G_i)$.

Now consider an HNN extension $G=G_1*_{G_e}$. Suppose that  an automorphism
equal to the identity on $G_1$ and sending the stable letter $t$ to $z_1tz_2\mi$
is a conjugation $i_m$.

Assume for a moment that $m\in G_1$.
We then write $u=z_2\mi m=t\mi z_1\mi
mt\in G_1\cap  t\mi G_1t=G_e$, so that $z_1=mt u\mi t\mi$ and
$z_2=mu\mi$   with $m\in Z(G_1)$ and $u\mi\in Z(G_e)$ (note that
$u$ commutes with $G_e$ because $m$ and $z_2$ do). This expresses that
$(z_1,z_2)$
belongs to   $J$.

Now we show $m\in G_1$. This is clear if  both inclusions
$G_e\hookrightarrow G_1$
are proper,   since $N_G(G_1)=G_1$.
If   an inclusion $G_e\hookrightarrow G_1$ is
onto, write $G=<G_1,t\mid tgt\mi=\varphi (g) \text{ if } g\in G_1>$, and
$m=t^{-k}gt^{\ell}$ with $g\in G_1$, $k \ge0$, and $k$ minimal. From $z_1
tz_2\mi=mtm\mi$ we get $g=t^kz_1(tz_2\mi t\mi)t^{-k}(tgt\mi)$, so that $g\in t
G_1t\mi$ if $k>0$. Minimality of $k$ then implies $k=0$. Writing that
$m=gt^{\ell}$ centralizes $G_1$ shows that $\varphi ^{|\ell|}$ is onto.
Since $\Gamma $ is not a
mapping torus, we have  $\ell=0$, and   $m\in G_1$.

{\it $\bullet$ Second step:
 $\Gamma $ has more than one edge and
there
is an edge $e$ such that each component of $\Gamma \setminus\{e\} $ is a minimal
graph of groups but not a mapping torus. }

First assume that
$e=v_1v_2$  separates. We write
$G=G_1*_{G_e}G_2$, where $G_i$ is the fundamental group of a   subgraph
$\Gamma _i$.    By minimality $G_e$ is a proper subgroup of $G_1$ and $G_2$, and
therefore $N_G(G _i)=G _i$.

A twist around an edge of $\Gamma _i$  may be represented
by an automorphism sending
$G _i$ to $G_i$ and equal to the identity on $G_{3-i}$. Similarly, a twist
 around $e$ near $v_i$ may be represented
by an automorphism sending
$G _i$ to $G_i$ by an inner automorphism and equal to the identity on
$G_{3-i}$.
With these representatives, consider the product $D_i$ of all twists by $z_f$
around edges $f$ contained in $\Gamma _i$. It sends
$G_i$ to itself, and because $\theta (\underline z)=1$ it coincides on
$G_i$ with
an inner automorphism of $G$.
Since $N_G(G_i)=G_i$, we deduce that $D_i$ restricts to an inner
automorphism of $G_i$.   Applying the induction hypothesis to $\Gamma _1$ and
$\Gamma _2$, we can then change $\underline z$ by an element of $J$ so as
to make every
$z_f$ trivial for
$f\inc \Gamma _1\cup\Gamma _2$ (this may change $z_e$, $z_{\bar e}$ by an
element
of the center of $G_{v_1}$ or $G_{v_2}$).

We   have now reduced to   the automorphism conjugating
each $G_i$ by some
$z_i\in Z_{G_{i}}(G_e)$. We have seen in the first step of the
proof that $z_i=mu_i$ with $m\in Z(G_e)$ and $u_i\in Z(G_i)$. The element $u_i$
fixes a point in the Bass-Serre tree of $\Gamma _i$ (it belongs to $G_{v_i}$).
Since $\Gamma _i$ is minimal,  $u_i$ acts trivially on the tree
and therefore belongs to  every edge and vertex group of $\Gamma _i$.
We conclude that $\underline z$ is the image by $j$ of the element of
$ \prod_{v\in V}Z(G_v)\times
 \prod_{f\in \E }
Z (G_f)$ whose component in $Z(G_v)$ is $u_i$ for $v\in \Gamma _i$, and whose
component in $Z(G_f)$ is $u_i\mi$ if $f\inc \Gamma _i$ and $m$ if $f=e$.

The argument when $e$ does not separate is similar, using representatives
sending $\pi _1(\Gamma \setminus\{e\})$ to itself (by the identity for twists
around $e$).   Minimality of
$\Gamma _1$ implies that both inclusions $G_e\hookrightarrow G_1$ are proper, so
$N_G(G_1)=G_1$.

The proof is now complete for graphs such that all edge groups are proper
subgroups of vertex groups. The general case requires one more argument.

{\it $\bullet $ Third step: graphs with a collapsible edge.}

We   prove the result for a graph $\Gamma $ containing a vertex $v$ of valence
$2$, with adjacent edges
$e,f$, such that $G_e=G_v$.  We shall make $\underline z$
trivial by multiplying it by elements of $J$ (i.e\. applying the edge and vertex
relations).

Call $u$ (resp\. $w$) the other vertex of $e$ (resp\.
$f$). Let $\Gamma '$ be the graph of groups obtained by collapsing $e$, i.e\.
replacing
$e\cup f$ by a single edge $g$ with group $G_f$.

The automorphism
$\theta (\underline z)$ twists by some $z_1$ around $e$ near $u$, by $z_2$ near
$v$, by $z_3$ around $f$ near $v$, by $z_4$ near $w$. This is the same
automorphism as twisting around $g$ in $\Gamma '$,  by $z_1z_2\mi z_3$ near
$u$ and
by
$z_4$ near
$w$ (and twisting as in $\Gamma $ around edges other than $g$). Since
$\Gamma '$ is
minimal, not a mapping torus, we may use the induction hypothesis and  kill
$z_1z_2\mi z_3$,
$z_4$, and all twists around edges other than $g$, by applying all edge and
vertex
relations of $\Gamma '$. If we use only those relations that exist in
$\Gamma $, we
cannot use the relation corresponding to the edge $g$ and we can only replace
$z_1z_2\mi z_3$ and
$z_4$ by some $m\in Z(G_f)$.

Going back to $\Gamma $, the relations let us replace
$(z_1, z_2,z_3,z_4)$ by
$(mz_3\mi z_2,z_2, z_3,m)$   and make all other $z$'s trivial.
 The element $mz_3\mi z_2$ commutes with $G_e$, and belongs to $G_e$
because all three factors do. It is then easy to make
$mz_3\mi z_2,z_2, z_3,m$ trivial by
  using the edge relations across $e$ and $f$ and the vertex relation at $v$.

{\it $\bullet $ Fourth step: conclusion.}

Let $\Gamma $ be any graph with more than one edge, satisfying the hypotheses
of the proposition. We may assume that there is no  vertex $v$ as in step
3.  There is then an edge $e$ with $\Gamma \setminus \{e\}$  minimal. We
argue as
in step 2.

If no component of
$\Gamma \setminus \{e\}$ is a mapping torus, we are done. If $e$ separates and
$G_i$ is a mapping torus, we observe that  the restriction of  $D_i$ to
$G_i$ is
 conjugation by an element  of
$G_{v_i}$ because all representatives used in step 2 send $G_1$ (resp\.
$G_2$) to
itself. Thus  step
$1$ applies (since mapping tori were ruled out only to prove $m\in G_1$).

There are now two cases left. First, $\Gamma $ may consist of two loops at a
vertex $v$, with all four edge groups mapping onto $G_v$. Then
$G=<G_v,t_1,t_2\mid
t_igt_i\mi=\varphi _i(g) \text{ if } g\in G_v>$. If an automorphism sending
$t_i$
to $z_it_iz_i'$ is $i_m$, then $m\in G_v$ because it normalizes $<G_v,t_i>$
and we
apply step 1. The last case is when $\Gamma $ has $3$ edges, each going from a
vertex
$v$ to a vertex $w$ (and all edge groups map onto vertex groups). We then use
vertex relations to kill twists at both ends of one of the edges, and we are
reduced to the previous case.
\cqfd\enddemo

\remark{Remark} Trivial examples show that the minimality hypothesis is
necessary. In the case of a mapping torus, the result fails precisely  in the
following situation: $\varphi \in\Aut(G_1)$ satisfies $\varphi ^\ell=i_g$
for some
$\ell\ge2$ and
$g\in G_1$ such that $g\mi\varphi (g)\in Z(G_1)$ but there is no $n\in
Z(G_1)$ with
  $g\mi\varphi (g)=n\mi\varphi (n)$. We do not know whether this situation is
possible.
\endremark

\head  {\sect Relative automorphism groups  }\endhead

  Many arguments in the beginning of this section are simple
computations which are left to the reader.

Let $H$ be a group, with a finite collection of (peripheral) subgroups
$\H=\{H_1,\dots,H_k\}$ (some groups may be repeated). For instance, $H$ is a
vertex group in  a graph of groups, and the $H_i$'s are the incident edge
groups.

Let
$\auta(H)$ be the  group of automorphisms which act on each $H_i$ as
conjugation
by some $a_i\in H$, and let $\outa(H)$ be its image in $\Out(H)$.

We
define $\autr(H)$ as the set of $(\alpha ;a_1,\dots,a_k)$, where $\alpha
\in\auta(H)$, $a_i\in H$,  and $\alpha (g)=a_iga_i\mi$ for $g\in H_i$. If
$\alpha $
acts on $H_i$ as conjugation by $a_i$, and $\beta $ acts as conjugation by
$b_i$,
then $\alpha \beta $ acts as conjugation by $\alpha (b_i)a_i$. We therefore
define
$$(\alpha ;a_1,\dots,a_k)(\beta ;b_1,\dots,b_k)=(\alpha \beta ;\alpha
(b_1)a_1,\dots,\alpha (b_k)a_k).
$$

This makes $\autr(H)$ into a group, with an exact sequence
$$\{1\}\to\prod_{i=1}^k Z_H(H_i)\to\autr(H)\to\auta(H)\to\{1\}.
$$

The map $g\mapsto(i_g;g,\dots,g)$ embeds $H$ as a normal subgroup of $\autr(H)$
and we define $\outr(H)$ as the quotient (note that $\outa(H)$ and $\outr(H)$
depend on $\H$). Up to isomorphism,
$\autr(H)$ and  $\outr(H)$ do not change if each $H_i$ is replaced by a
conjugate.
One easily checks:

\nom\suitex
\thm{Lemma \sta} There is  an exact sequence
$$\{1\}\to\biggl(\prod_{i=1}^k
Z_H(H_i)\biggr)/Z(H)\to\outr(H)\overset {\pi _H}\to\to\outa(H)\to\{1\},
$$
with $Z(H)$ embedded diagonally into $\prod Z_H(H_i)$. If $Z_H(H_i)\inc H_i$
for every $i$, this extension is central. \vskip-10pt \cqfd
\fthm

\example{ Examples}

$\bullet$ If $k=0$, then $\outr(H)=\Out(H)$. If $k=1$, the map $(\alpha
;a)\mapsto i_a\mi
\alpha $ identifies
$\outr(H)$ with the group of automorphisms of $H$ equal to the identity on
$H_1$.

$\bullet$ Suppose the $H_i$'s are all trivial. Then $\outa(H)=\Out(H)$. If
$k=1$,
we have
$\outr(H)=\Aut(H)$ and the  exact sequence of Lemma \suitex{} is simply
$\{1\}\to H/Z(H)
\to\Aut(H)\to\Out(H)\to\{1\}$. If $k\ge1$, we get
$\{1\}\to H^k/Z(H)
\to\outr(H)\to\Out(H)\to\{1\}$, with $Z(H)$ embedded diagonally into $H^k$.

$\bullet$  If $H=\Z$ and at least one $H_i$ is nontrivial, then
$\outr(H)\simeq\Z^{k-1}$, the quotient of $\Z^k$ by $\Z$ embedded diagonally.

$\bullet$ If $H$ is the fundamental group of a
compact  surface $\Sigma $ with $  \chi (\Sigma
)<0$, and
$\H$ comes from the boundary curves, then $\outr(H)=\outr(\Sigma )$, the group
of
 homeomorphisms fixing $\bo\Sigma $ pointwise, up to isotopy   relative
to
$\bo\Sigma $. It is a central extension of $\outa(\Sigma )=\outa(H)$
by
$\Z^k$, where $k$ is the number of boundary components.

In general this extension is not trivial. For instance, the ``lantern relation''
(see e.g\. [\IMc, section 4]) implies that the kernel $\Z^k$ is contained in the
derived subgroup of
$\outr(\Sigma )$ whenever $\Sigma $ has genus $\ge2$.  When $\Sigma $
is a once-punctured torus the extension, though not trivial, is trivial above a
finite index subgroup of $\outa(\Sigma )$. This is probably not true for
surfaces
of high genus.
\endexample

Returning to a graph of groups $\Gamma $ as in section 2, we consider
$\autr(G_v)$
and $\outr(G_v)$ relative to the collection of incident edge groups ($k$ is the
valence of $v$ in $\Gamma $). We denote by
$\ds\pi_v:
\outr(G_v)\to
  \outa(G_v)$ the natural projection, and we define
$ \pi=\prod_{v\in V} \pi _v: \prod
\outr(G_v)\to
\prod  \outa(G_v)$.

The extension construction of section 2  may be viewed as a
 homomorphism
$\autr(G_v)\to\Aut(G)$. It induces
$\lambda _v:\outr(G_v)\to\Out^\Gamma _0(G)$, with  $\rho _v\circ \lambda _v=\pi
_v$. The
images of
$\lambda _v$ and $\lambda _w$ commute for $v\neq w$, so that we get
$\ds \lambda :\prod_{v\in V} \outr(G_v)\to \Out^\Gamma _0(G) .$
We call $\Out^\Gamma _1(G)$ the image of $\lambda $.

\nom\suite
\thm{Proposition \sta}
There is a commutative triangle of exact sequences
$$
\xymatrix
{N\ \ar@{^{(}->} [dr]
\\ \prod \frac{\prod_{e\in E_v }
Z_{G_v}(G_e)}{Z(G_v)}\ \ar @{^{(}->}[r]&{\prod\rmcg(G_v)}\ar @{>>}[rr]^{\pi
}\ar@{>>} [rd]
_{\lambda }
&&{\prod\mcg(G_v)}\\ &&{\Out_1^\Gamma (G)}\ar@{>>} [ur]_{\rho_1 }\\
&\T\ \ar@{^{(}->} [ur]
}
$$
where $\Out^\Gamma _1(G)$ is a subgroup of
$\Out^\Gamma  (G)$,  all products are taken over $  V$, and $\T$ is the group of
twists.

 If all
groups
$\Out(G_e)$ are finite, then $\Out^\Gamma _1(G)$ has finite index in
$\Out^\Gamma
(G)$. If all edge groups are trivial, then $\Out^\Gamma _1(G) =\Out^\Gamma _0
(G)$ and $\lambda $ is an isomorphism.
\fthm

\demo{Proof} We have seen that $\pi $ factors as $\rho_1\circ \lambda $,
with $\rho
_1$ the restriction of $\rho $ to $\Out^\Gamma _1(G)$. The kernel of $\pi $ is
given by Lemma \suitex{}, and  the extension construction shows
$\ker\rho _1=\T$.  Note that $N$ is the kernel of the natural epimorphism
  $\ds\kappa :\prod_{v\in V}\frac{\prod_{e\in E_v }
Z_{G_v}(G_e)}{Z(G_v)}\to
 \T$.

The group $\Out^\Gamma _1
(G)$ is mapped by $\rho $ onto $\prod\outa(G_v)$, and it contains
$\T$. If all groups
$\Out(G_e)$ are finite,   Proposition \indice{} implies that $\Out^\Gamma _1
(G)$ has finite index in $\Out^\Gamma_0
(G)$, hence in $\Out^\Gamma
(G)$. If every $G_e$ is trivial, then $ \outa(G_v)=
\Out(G_v)$, and $\T =\ker\rho $ by  Proposition \noyau{}. Thus $\Out^\Gamma _1
(G)=\Out^\Gamma_0
(G)$.
Furthermore $\lambda $ is injective,  because $N=\ker\lambda =\ker\kappa $ is
generated by centers of edge groups  by Lemma \relat{}.
\cqfd\enddemo

\head  {\sect Hyperbolic groups  }\endhead

In this section $G$ is a one-ended hyperbolic group. The {\it elementary
subgroups\/} of $G$ are the finite or virtually cyclic (= $2$-ended) subgroups.
Recall (see [\MNS]) that there are two types of infinite elementary
groups, those
with infinite center (which map onto
$\Z$) and those with finite center (which map onto the infinite dihedral group).
The normalizer of a virtually cyclic subgroup is a maximal elementary
subgroup. If
$G $ is torsion-free, nontrivial elementary subgroups are cyclic.

 Let  $\Gamma $ be the
JSJ splitting of $G$, as constructed by Bowditch [\Bow].
Edge groups of $\Gamma $ are virtually cyclic. There are  three types of
vertices:

1.  {\it Elementary vertices\/}, whose group is a maximal elementary subgroup.

2. {\it Orbifold vertices\/}, whose group is a
  ``hanging Fuchsian group'', the fundamental group of  a $2$-orbifold.

3. {\it Rigid vertices\/}, whose group $G_v$ cannot be split any further over an
elementary subgroup. By Paulin's theorem [\Pau] and Rips theory [\BF], the group
$\outa(G_v)$ is finite (we need a relative version of these results, which
follows from Theorem 9.6 of [\BF] since applying [\Pau] to automorphisms
of $\outa(G_v)$ gives a tree in which every peripheral subgroup of
$G_v$ fixes a point).

We denote by $V_1$, $V_2$, $V_3$ the corresponding sets of vertices of
$\Gamma $.
Every edge joins an elementary vertex to an orbifold or rigid   vertex, we
denote
$\E=\E_2\cup\E_3$ the   set of geometric edges. We shall add the symbol
$\infty$ to
indicate that we restrict to  edges or elementary vertices whose group has
infinite center (orbifold and rigid vertex groups always have finite
center). The
notation
$|A|$ denotes the number of elements  of a finite set
$A$.

Edge groups may fail to be maximal as elementary subgroups of $G$, but it is
easy to see that they are maximal elementary in the corresponding rigid or
orbifold
vertex group (otherwise their normalizer would not be elementary). If $v\in
V_2\cup V_3$, we have
$Z_{G_v}(G_e)=Z(G_e)$ for $e\in E_v$, and  $\outr(G_v)$ is a central
extension of
$\outa(G_v)$ by Lemma \suitex.

\nom\exun
\thm{Theorem \sta} There is an exact sequence
$$1\to\T\to\Out  _2(G)
\overset {\rho _2}\to
\to \prod_{v\in V_2}
\outa(G_v)\to1,
$$
where $\Out  _2(G)$ has finite index in $\Out  (G)$.
The group of twists $\T$  is virtually $\Z^n$, with
$n=|\E^\infty|-|V_1^\infty|$. If $G$ is torsion-free, then $\T$ is free
abelian of
rank $|\E|-|V_1|$ and the extension is central.
\fthm

\demo{Proof}    Since
$\outa (G_v)$ is finite if $v\in V_1\cup V_3$, we
define
$ \Out
_2(G)$ as the kernel of the map
$\Out  _1(G)\to
\prod_{v\in V_1\cup V_3 } \outa(G_v)$.
The exact sequence then comes from Proposition \suite.

We use Proposition
\relat{} to compute the rank of $\T$.   Since  $Z_{G_v}(G_e)=Z(G_e)$ if
$e$ is an
edge and
$v$ is its rigid or orbifold vertex, we see that $\T$ is the quotient of
$\ds\T'=\prod_{v\in
V_1}
\frac{\prod_{e\in E_v}Z_{G_v}(G_e)}{Z(G_v)}$ by a finite group (corresponding to
orbifold and rigid vertex relations).

First suppose that $G$ is
torsion-free. Then $\T=\T'$ because $Z(G_v)$ is trivial if $v\notin V_1$. Of
course $Z_{G_v}(G_e)=Z(G_v)=G_v=\Z$ for $v\in V_1$, and the result follows,
recalling that
$Z(G_v)$ is embedded diagonally in
$\prod Z_{G_v}(G_e)$.

If $G$ has torsion, we obtain that $\T'$, hence also $\T$, is virtually $\Z^n$
with $n=|\E^\infty|-|V_1^\infty|$
because
$Z_{G_v}(G_e)$ is infinite if and only if $G_e$ has infinite center.

The extension is central as soon as  $\outr(G_v)$ is a central extension of
$\outa(G_v)$ for every $v $, in particular if $G$ is torsion-free.
\cqfd\enddemo

   Proposition \suite{} also yields the exact sequence
$$\{1\}\to N
\to A\times \prod_{v\in V_2}
\outr(G_v)
\to\Out  _2(G)\to\{1\},
$$
where $\ds A=\prod_{v\notin V_2}
\ker\pi _v=
\prod_{v\notin V_2}\frac{\prod_{e\in E_v }
Z_{G_v}(G_e)}{Z(G_v)}$ is virtually abelian. Note that $N$ is also virtually
abelian (because it is a subgroup of $\ker\pi $). We shall now use the
edge and vertex relations in order to replace
$A$ and
$N$ by simpler groups.

We define a set $Q$ of edges as follows. We consider   $\E_3^\infty$ (edges $e$
connected to a rigid vertex, with $Z(G_e)$   infinite), and for each
$v\in V_1$ whose group has infinite center we discard one   edge of
$\E_3^\infty$ adjacent to $v$ (if there is at least one). We view $Q$ as a
set of oriented edges, going   from $V_1$ to $V_3$. Note that
$q=|Q|$ equals
$|\E_3^\infty|-|V_1^\infty|+r$, where $r$ is the number of vertices of
$|V_1^\infty|$ connected only to orbifold vertices.

Now consider the map $$\psi :\prod_{e\in Q} Z_{G_{o(e)}}(G_e)\times \prod_{v\in
V_2}\outr(G_v)
\to\Out  _2(G)$$
obtained by mapping the first factor to $A$ and applying $\lambda $.

\nom\disj
\thm{Lemma \sta} The image of $\psi $ has finite index. The images of the two
factors have finite intersection. The  intersection of $\ker\psi $
with
the first (resp\. second) factor is finite (resp\. virtually $\Z^r$).
\fthm

\demo{Proof} The first factor maps into $\T$, the kernel of $\rho _2:\Out  _2(G)
\to \prod_{v\in V_2}
\outa(G_v)$. The image
by
$\psi
$ of the second factor is mapped onto the whole of $\prod_{v\in V_2}
\outa(G_v)$ by $\rho _2$, and
the intersection of $\psi \mi(\T)$ with the second factor is the subgroup
$\ds\prod_{v\in V_2}\frac{\prod_{e\in E_v } Z_{G_v}(G_e)}{Z(G_v)}$. It therefore
suffices to prove the statements of the lemma for the map
$$\prod_{e\in Q} Z_{G_{o(e)}}(G_e)\times \prod_{v\in V_2}\frac{\prod_{e\in E_v }
Z_{G_v}(G_e)}{Z(G_v)}\to\T,$$
or   for
$$\prod_{e\in Q} Z_{G_{o(e)}}(G_e)\times \prod_{v\in V_2} {\prod_{e\in E_v }
Z_{G_v}(G_e)} \to\T$$ since $Z(G_v)$ is finite for $v\in V_2$. But this follows
directly from Proposition
\relat, as we now explain.

In this argument, we neglect finite groups (and we don't distinguish between a
group and a finite index subgroup).  We know that $\T$ is generated by all
groups $Z_{G_{o(e)}}(G_e)$. The edge relations reduce the generating set
to
$\prod_{e\in \bar Q} Z_{G_{o(e)}}(G_e) \times  \prod_{v\in V_2}
{\prod_{e\in E_v }
Z_{G_v}(G_e)}$, where $\bar Q$ is the set of oriented edges from $V_1$ to $V_3$
whose group has infinite center. Now we have to consider vertex relations at
vertices $v\in V_1^\infty$. Vertices connected to at least one rigid vertex
let us reduce from $\bar Q$ to $Q$, while vertices connected only to
orbifold vertices generate the desired $\Z^r$ in the second factor.
\cqfd\enddemo

\thm{Theorem \sta} Let $r$ be the number of vertices of
$|V_1^\infty|$ connected only to orbifold vertices. Let $q=
|\E_3^\infty|-|V_1^\infty|+r$ and $s=|\E_2^\infty|-r=n-q$.
The group $\Out (G)$ is virtually a direct product $\Z^q\times M$, where $M$
fits in exact sequences
$$\{1\}\to Z_r
\to  \prod_{v\in V_2}
\outr(G_v)
\to M\to\{1\}
$$
$$\{1\}\to Z_s
\to M\to \prod_{v\in V_2}
\outa(G_v)
\to\{1\}
$$
with $Z_r$ (resp\. $Z_s$)   virtually $\Z^r$ (resp\. $\Z^s$). If $G$ is torsion
free, then $Z_r=\Z^r$, $Z_s=\Z^s$, and the extensions are central.
\fthm

\demo{Proof}
This follows from Lemma \disj, with
 $M$  the image by $\psi $ of the second factor, and $\Z^q$ the image of a
$\Z^q$ of finite index in $\prod_{e\in Q} Z_{G_{o(e)}}(G_e)$. The group $Z_s$ is
virtually $\Z^s$ because $\Z^q\times Z_s$ has finite index in $\T$, which is
virtually $\Z^{q+s}$.

If $G$ is torsion-free, then $Z_r$ is contained in the
free abelian group $\ds\prod_{v\in V_2}\frac{\prod_{e\in E_v }
Z_{G_v}(G_e)}{Z(G_v)}$ and $Z_s$ is contained in $\T=\Z^n$. The extensions
are central by Lemma \suitex{} and Theorem \exun.
\cqfd\enddemo

\thm{Corollary \sta} Suppose $G$ is torsion free.
   If every
cyclic vertex   is connected to at least one rigid vertex, then $\Out(G)$ is
virtually a direct product of a free abelian group with mapping class
groups $\outr(\Sigma _v)$ of surfaces with boundary. \cqfd
 \fthm

\head  {\sect Hyperbolic groups with $\Out(G)$ infinite }\endhead

In this section we prove Theorem \conj. The following fact  is quite general
(compare [\MNS, Proposition 2.1]).

\thm{Proposition \sta} Let $G$ be any group. Suppose $G=G_1*_HG_2$ or
$G=G_1*_H$, where $H$ has infinite center but $G_1$, $G_2$ have finite center.
Then $\Out(G)$ is infinite.
\fthm

 This may be viewed as a special case of Proposition \relat{} (see also
[\MNS] for
the case of an amalgam), and proves one direction of Theorem \conj{} (note that
$G_1$ has finite center if $G=G_1*_H$ is hyperbolic and $H$ is virtually
cyclic).

Conversely, suppose $G$ is one-ended, hyperbolic, with $\Out(G)$ infinite.
Consider the exact sequence of Theorem \exun. Since $G$ splits as in Theorem
\conj{} if $n>0$, we may assume that there exists $v\in V_2$ with $\outa(G_v)$
infinite (note that this argument requires Proposition   \indice, but
not Proposition \relat).

The group $G_v$ maps into the group of isometries of $\bold H^2$ with finite
kernel (see [\Bow]). We start by assuming this kernel to be
trivial.

The (convex core of the) quotient orbifold
$\Sigma _v=\bold H^2/G_v$ is homeomorphic to a compact surface, with a singular
set which may be
one-dimensional: this phenomenon is responsible for the existence of   virtually
cyclic subgroups of
$G_v$ with finite center. The required splitting of $G$ will come from a
splitting
of $G_v$, given by an essential $2$-sided simple closed curve $C\inc \Sigma
_v$ disjoint from the singular set.

If $G_v$ is torsion-free, then $\Sigma _v$ is a regular compact hyperbolic
surface
with infinite mapping class group. Thus $\Sigma _v$ cannot be a pair of pants or
a twice-punctured projective plane.
On all other surfaces we can find the required essential $2$-sided curve $C$
(essential means that, if
$C$ separates, then both  complementary regions have
negative Euler characteristic).

If $G_v$ has torsion, then $\Sigma _v$ is an orbifold. By [\HM], generalized
to the nonorientable case in [\Fuj], the group $\outa(G_v)$ is
commensurable with
the mapping class group of a regular compact surface $\Sigma _0$ obtained from
$\Sigma _v$ by removing an open neighborhood of the singular set. An essential
curve $C\inc
\Sigma _0$ yields the desired splitting of $G_v$.

Up to now we have assumed   that $G_v$ acts effectively on $\bold H^2$. If not,
there is a finite normal subgroup $F\inc G_v$ with $G'_v=G_v/F
\inc Isom(\bold H^2)$. As pointed out in [\Bow], $F $ is the unique maximal
finite normal subgroup of $G_v$. Automorphisms of $G_v$ therefore induce
automorphisms of $G'_v$. The map $\Aut( G_v)\to\Aut(G'_v)$, hence also the map
$\outa( G_v)\to\outa(G'_v)$, is at most $|F|^k$-to-one if $G_v$ may be
generated by
$k$ elements, and therefore  $\outa(G'_v)$ is infinite.  The splitting of
$G'_v$ constructed above lifts to a
  splitting of $G_v$ and extends to a splitting of $G$ (the splitting of
$G_v$ is
over a group that maps onto $\Z$, hence has infinite center).

This completes the proof of Theorem \conj.

\remark{Remark} We sketch a proof which does not use the orbifold theory of
 [\HM] and [\Fuj]. View $\Sigma _v$ as the quotient of a surface $\Sigma $ by a
finite group $\Omega $. An algebraic argument shows that the mapping class
group of $ \Sigma  $ contains an element of infinite order commuting with
$\Omega $.  Represent it by a homeomorphism consisting of  pseudo-Anosov
homeomorphisms of subsurfaces and Dehn twists along disjoint curves. Because of
$\Omega $-equivariance,  the support of these maps is disjoint from the
one-dimensional components of the branching locus of $\Omega $, and we can
construct a curve
$C$.
\endremark

\head  {\sect Automorphisms of finite order }\endhead

This section is devoted to the proof of  Theorem \nombfi.

\subhead Algebraic preliminaries \endsubhead

We denote by $\cc$
the class of groups with only finitely many conjugacy classes of torsion
elements.

\nom\algeb
\thm{Lemma \sta} Let $ B$ be a group. If  one of the following holds, then
$B\in\cc$.
\roster
\item $B=A\times A'$, or $B=A*A'$, with $A,A'\in\cc$.
\item   $B$ is a subgroup of finite index of $A$, and $A\in\cc$.
\item There is an epimorphism $p:B\to A$, with $A\in\cc$, such that
$p\mi(F)\in\cc$ for every finite subgroup $F\inc A$.
\item $B$ contains a normal subgroup $H$ with $H$ hyperbolic and $B/H\in\cc$.

\endroster
\fthm

\demo {Proof} (1), (2), (3) are simple  observations. By (3), it suffices
to prove
(4) when $B/H$ is finite. But then $B$ is hyperbolic, hence in $\cc$.
\cqfd\enddemo

\example{Remark} Unfortunately, there is no converse to (2). Let $A$ be the
wreath
product $\Z\wr \Z$, with commuting generators $a_i$ ($i\in\Z$)  and an extra
generator
$t$ conjugating $a_i$ to $a_{i+1}$. Let $B$ be the semi-direct product
$A\rtimes_\varphi   \Z_2 $, where $\varphi $ is the involution sending $a_i$ to
$a_i\mi$ and $t$ to $t$. Then $B$ is a finitely generated group with
infinitely
many conjugacy classes of elements of order
$2$, but $A$ is a torsion-free subgroup of index $2$.
\endexample

\nom\algebb
\thm{Proposition \sta} Let $ B$ be a group. If  one of the following holds (with
$n$ a positive integer), then $B\in\cc$.
\roster
\item "{(i)}" There is an exact sequence $\{1\}\to K
\to B\to A
\to\{1\}
$, with $K$ virtually $\Z^n$ and $A\in\cc$.
\item "{(ii)}" $B$ is the semi-direct product $A^n\rtimes \s_n$, with the
symmetric
group
$\s_n$ acting on $A^n$  by permuting the factors, and $A\in\cc$.
\endroster
\fthm

\demo{Proof} We first prove (i)  when
$K=\Z^n$ and $A$ is finite. For every $a\in A$, choose a lift $b_a\in B$ of
finite
order, if there is one. Also fix an integer $q$ such
that $b^q=1$ for every torsion element $b\in B$.

Given $b$ of finite order, we write $b=tb_a $ with $a\in A$
and $t\in\Z^n$, and we
consider the relation
$(tb_a)^q=1$.
The action of $b_a$ on $\Z^n$ by conjugation is given by a matrix $P\in
GL(n,\Z)$,
with $P^q=1$, and $(tb_a)^q=tP(t)P^2(t)\dots P^{q-1}(t)b_a^q$.  We get
$tP(t)P^2(t)\dots P^{q-1}(t)=1$.

Over the rationals we
would conclude that $t $ may be written $s\mi P(s)$ for some $s $, because
the restriction of
$P$ to the subspace generated by   $t,P(t),P^2(t),\dots,P^{q-1}(t)$
does not have $1$ as an eigenvalue. Here we    obtain $t=s\mi P(s) t'$ with $t'$
belonging to some finite set (depending on $a$). We then write
$b=tb_a=s\mi P(s )t'b_a=s\mi t'P(s )b_a=s\mi t'b_as $. Since there are finitely
many couples $(b_a,t')$, we get $B\in\cc$.

We have now proved that every group which is virtually $\Z^n$ belongs to $\cc$.
The general case of (i) follows from   assertion (3) of Lemma
\algeb{}.

To prove (ii), write $b\in B$ as $(a_1,\dots,a_n)\sigma $ with $a_i\in A$ and
$\sigma \in\s_n$.
  First suppose that
$\sigma $ is an
$n$-cycle, say the standard $n$-cycle. For $c_i\in A$ we have
$\sigma
(c_1\mi,\dots,c_n\mi)=(c_2\mi,\dots,c_n\mi,c_1\mi)\sigma$ and
$$(c_1,\dots,c_n)(a_1,\dots,a_n)\sigma
(c_1,\dots,c_n)\mi=(c_1a_1c_2\mi,\dots,c_{n-1}a_{n-1}c_n\mi,c_na_nc_1\mi)\sigma.
$$ It follows that $b$ is conjugate to $(a ,1,\dots,1)\sigma $, with $a=a_1\dots
a_n$, and that $(a ,1,\dots,1)\sigma $ and $(a' ,1,\dots,1)\sigma $ are
conjugate
whenever $a$ and $a'$ are conjugate in $A$. Also note that
$[(a ,1,\dots,1)\sigma ]^n=(a,\dots,a)$, and therefore
$(a ,1,\dots,1)\sigma $ has finite order if and only if $a$ has finite order. We
deduce that $B$ contains only finitely many classes of torsion elements whose
associated permutation is an $n$-cycle.

The general case follows easily. Write $\sigma =\sigma _1\dots\sigma _k$ as a
product of cycles with disjoint support, including cycles of length $1$. This
induces a   decomposition of $b$ as a product of commuting factors $b=(\theta
_1\sigma _1)\dots(\theta _k\sigma _k)$, where the $j$-th component of
$\theta _i\in
A^n$ is $a_j$ if $j$   belongs to the support of $\sigma _i$ and $1$
otherwise. We
conclude by arguing as above in each factor separately.
\cqfd\enddemo

\subhead Groups with one end \endsubhead

Let $G$ be a  torsion-free  one-ended hyperbolic group, and   $\Gamma
$ its JSJ splitting   as before. We associate to
$\hat\alpha
\in\Out(G)$ the induced permutation $\sigma $ on the vertex set $V$ of $\Gamma
$. This defines  $\mu :\Out(G)\to\s(V)$, and we  call $S$ its image.

In section 2
we restricted to  automorphisms $\alpha $ acting trivially on
$V$, and we defined
$\rho _v(\hat\alpha )\in\Out(G_v)$ for $v\in V$. In the general case,  we
obtain  a
set of
  isomorphisms $  \eta _v(\hat\alpha ):G_v\to G_{\sigma (v)}$, each defined
up to composition with inner automorphisms. These isomorphisms preserve the
peripheral structure (the set of   conjugacy classes of edge groups), but in
general not every isomorphism preserving the peripheral structure arises in this
way.

 We view this
construction as a homomorphism $\eta $ from $\Out(G)$ to a group $\Lambda
$  defined as follows.
  An element of $\Lambda $
consists of
$\sigma
\in S$ and isomorphisms $\eta _v:G_v\to G_{\sigma (v)}$ preserving the
peripheral
structure, defined up to inner automorphisms.  The group law on $\Lambda $ is
defined in the obvious way.

The rest of the proof consists in showing that {\it $\Lambda \in\cc$, the
image  of
$\eta
$ has finite index in
$\Lambda $, and $\ker\eta $ is virtually $\Z^n$.}
These facts imply Theorem \nombfi{} (for $G$ one-ended) by \algeb.2 and
\algebb.i.

The kernel of the natural projection $\pi :\Lambda \to S$ is $\prod_{v\in
V}A_v$,
with
$A_v\inc \Out(G_v)$ consisisting of automorphisms preserving the peripheral
structure.
The group $A_v$ is finite
if $v$ is a rigid or cyclic vertex. If $v$ is a surface vertex, then $A_v$ is
the full mapping class group of the surface (boundary components may be
permuted,  orientation may be reversed). In all cases $A_v\in\cc$.

The reason for introducing the abstract group $\Lambda $ is that
the map $\pi $ always has a
section: for every orbit  $V_i$  of the action of $S$ on $V$, choose
isomorphisms
(compatible with the peripheral structures) of the groups $G_v$, $v\in
V_i$, with a
fixed model, and use these identifications to construct the section.

We now know that $\Lambda $ is a semidirect product $(\,\prod_{v\in
V}A_v)\rtimes
S$. It embeds as a finite index subgroup into $ \prod _i\Lambda _i$, with
$\Lambda _i=(\,\prod_{v\in V_i}A_v)\rtimes \s(V_i)$.
Proposition \algebb.ii{} implies $\Lambda _i\in\cc$, and
$\Lambda \in\cc$ by \algeb.1 and \algeb.2.

Now consider   $\eta :\Out(G)\to \Lambda $.   By
Proposition \debut{} it contains $\prod_{v\in V}\mcg(G_v)$, which has
finite index
in $\prod_{v\in V}A_v$. Thus the image of $\eta $ has finite index. Its
kernel  contains $\ker\rho $ with finite index (an element of $\ker\eta
$ belongs to $\ker\rho
$ if and only if it acts trivially on the set of edges of $\Gamma $).
But  $\ker\rho $ is virtually $\Z^n$ by Proposition \indice{} and Theorem \exun.

This completes the proof of Theorem \nombfi{}  for one-ended groups.

 \subhead Groups with infinitely many ends \endsubhead

We now write
$G=G_1*\dots*G_p*F_q$ with $G_i$ one-ended and $F_q$ free of rank $q$. We
argue by
induction on the Kurosh rank $r(G)=p+q$.

The relative train track maps of Bestvina-Handel [\BH] have been
generalized to free products in [\CT], so that we can represent any $\hat\alpha
\in\Out(G)$   by a   relative train track map
$\varphi :X\to X$. We assume that $\hat\alpha $ has finite order, and we use
$\varphi
$ to construct an
$\alpha
$-invariant splitting of $G$.

 The map $\varphi $ is a self-homotopy equivalence of a
complex $X$ obtained by attaching connected complexes $X_i$ with $\pi
_1(X_i)\simeq
G_i$ onto vertices of a  finite connected graph $\Delta  $ with $\pi _1(\Delta)
\simeq F_q$.  Let $\Delta _t $ (resp\. $\Delta _t'$) be the union of edges
(resp\. open edges) in the top stratum. The complement of
$\Delta _t'$ in
$X$ is mapped to itself by
$\varphi $.

Since
$\hat\alpha
$ has finite order, $\Delta _t$ cannot be exponentially growing. It follows that
the (unoriented) edges of $\Delta _t$ are cyclically permuted (modulo the lower
strata):  they may be numbered
$e_1,\dots, e_k$ in such a way that $\varphi (e_i)$ consists of $e_{i+1}$,
possibly preceded and followed by a path disjoint from $\Delta _t'$.

To obtain an
$\alpha $-invariant graph of groups $\Omega _\alpha $, we first make
$\Delta $ into
a graph of groups  in the obvious way (vertex groups are
free products of $G_i$'s), and then we collapse all edges not in $\Delta _t$.

Edge groups of $\Omega _\alpha $ are trivial. Vertex groups $G_v$ are free
products
of
$G_i$'s and $\Z$'s with Kurosh rank $r(G_v)<r(G)$, so $\Out(G_v)\in\cc$ by
induction. Furthermore, there are only finitely many possible $\Omega
_\alpha $'s,
up to isomorphism. This means that there exist finitely many graph of groups
decompositions $\Omega _i$ of $G$, such that every torsion element of
$\Out(G)$ is
conjugate to an element of an $\Out^{\Omega _i}(G)$. It therefore suffices
to show:

\thm{Proposition \sta} Suppose that $G$ is the fundamental group of a
finite graph
of groups
$\Omega
$ with trivial edge groups. Assume that  vertex groups are torsion-free
hyperbolic
groups $G_v$ with $\Out(G_v)\in\cc$. Then $\Out^\Omega (G)  \in\cc$.
\fthm

\demo{Proof} As usual, we denote by $\Out^\Omega _0 (G)$ the kernel of the
natural
map
$\psi :\Out^\Omega  (G)\to\Aut(\Omega )$. Let $L$ denote the image. Let $V$
be the
vertex set of $\Omega $.

By Proposition \suite{}, we have $\Out^\Omega _0 (G)=\prod_{v\in
V}\outr(G_v)$. Here $\outr(G_v)$ is the automorphism
group of the vertex group $G_v$ relative to the family $\H$  consisting of the
trivial group repeated $n_v$ times ($n_v=|E_v|$ being the valence of $v$ in
$\Omega
$). It fits in an exact sequence
$$\{1\}\to (G_v)^{n_v}/Z(G_v)
\to\outr(G_v)\to\Out(G_v)\to\{1\}$$ (see the second example in section 4).
If $G_v=\Z$, then $\outr(G_v)$ is virtually $\Z^{n_v-1}$.
If $G_v\neq\Z$, the center of $G_v$ is trivial and the kernel is simply
$(G_v)^{n_v}$; the natural action of $\s(E_v)$ on $\outr(G_v)$ permutes the
factors
of $(G_v)^{n_v}$.

Let $V_i$ be the orbits of $V$ under the action of $L$. For each $v\in V_i$, we
choose an identification of $G_v$ with a fixed group $G_i$, and also a
numbering of
edges adjacent to $v$ (i.e\. a bijection $E_v\to\{1,\dots,n_v\}$). This gives
canonical isomorphisms $\outr(G_v)\to\outr(G_w)$ for $v,w$ in the same
$V_i$, and
provides a section to the map $\psi :\Out^\Omega  (G)\to  L$.

We now have an action of $L$ on $\prod_{v\in
V}\outr(G_v)=\ker\psi $, which is a restriction of the natural action of the
semi-direct product $\bigl(\prod _{v\in V}\s(E_v)\bigr)\rtimes \prod_i \s(V_i)$. In
other words,
$\Out^\Omega  (G)$ is a finite index subgroup of $\biggl(\prod_{v\in
V}\outr(G_v)\biggr)\rtimes\biggl(\bigl(\prod _{v\in V}\s(E_v)\bigr)\rtimes
\prod_i
\s(V_i)\biggr)$. It suffices to show that this larger group is in $\cc$.

We may rewrite it as $\ds\prod_i\biggl( \prod_{v\in V_i}\bigl( \outr(G_v)\rtimes
\s(E_v)\bigr)\rtimes \s(V_i)\biggr)$, and by \algeb.1 and \algebb.ii we
need only
show
$\outr(G_v)\rtimes \s(E_v)\in\cc$. We may assume $G_v\neq\Z$, since
otherwise the
group is virtually $\Z^{n_v-1}$, hence in $\cc$ by \algebb.i.

With the notations of section 4, the map $(\alpha
;a_1,\dots,a_k)\mapsto\bigl((\alpha ;a_1),\dots,(\alpha ;a_k)\bigr)$ induces
$\zeta :\outr(G_v)\to\prod_{i=1}^{n_v}\Aut(G_v)$ (recall that $\outr =\Aut $
when $\H$ consists of one copy of the trivial group). Since $G_v$ has trivial
center, $\zeta $ is an embedding. In general, its image does not have
finite index.

Now let $F$ be a finite subgroup of $\Out(G_v)$, and
$\outr_F(G_v)$, $\Aut_F(G_v)$ its preimages. The images of the embeddings
$\outr_F (G_v)\to\prod_{i=1}^{n_v}\Aut_F(G_v)$ and
$\outr_F (G_v)\rtimes
\s(E_v)\to\bigl(\prod_{i=1}^{n_v}\Aut_F(G_v)\bigr)\rtimes
\s(E_v)$ have finite index (all groups are virtually $(G_v){}^{n_v}$). The group
$\Aut_F(G_v)$ is hyperbolic, hence   in $\cc$. By \algeb.2 and \algebb.ii, we
deduce that
$\outr_F(G_v)\rtimes \s(E_v)\in\cc$ for every finite  $F\inc\Out(G_v)$. Applying
\algeb.3 to the map from $\outr(G_v)\rtimes \s(E_v)$ to $\Out(G_v)$, we get the
required result
$\outr(G_v)\rtimes
\s(E_v)\in\cc$.
\cqfd\enddemo

\bigskip

\Refs
\widestnumber\no{99}
\refno=0

\bref  \by
H. Bass, R. Jiang  \paper Automorphism groups of tree actions and of
graphs of groups\jour J. Pure Appl. Algebra \vol 112 \yr1996\pages 109--155
\endref

\bref \by M. Bestvina, M. Feighn \paper   Stable actions of groups
on real trees \jour Invent. Math. \vol121 \yr1995\pages 287--321
\endref

\bref \by M. Bestvina, M. Handel\paper Train tracks for automorphisms
of the free group \jour Ann. Math.\vol135 \yr1992\pages1--51 \endref

\bref \by B. Bowditch \paper Cut points and canonical splittings of
hyperbolic groups\jour Acta Math. \vol180\yr1998\pages145--186\endref

\bref\by D.J. Collins, E.C. Turner\paper Efficient representatives for
automorphisms of free products\jour Michigan Math. Jour.\vol41\yr1994\pages
443--464\endref

\bref \by K. Fujiwara \paper On the outer automorphism group of a hyperbolic
group\jour Israel Jour. Math. (to appear)
\endref

\bref \by W.J. Harvey,
     C. Maclachlan   \paper On mapping-class groups and
Teichm\"uller spaces\jour Proc. London Math. Soc. \vol 30 \yr1975\pages
496--512\endref

\bref\by N.V. Ivanov, J.D. McCarthy\paper On injective homomorphisms
between Teichm\"uller modular groups I\jour Invent. Math. \vol 135 \yr1999\pages
425--486
\endref

\bref\by G. Levitt, M. Lustig \jour in preparation
\endref

\bref \by C.F. Miller III, W.D. Neumann, G.A. Swarup \paper Some examples of
hyperbolic groups, {\rm pp\. 195--202 in ``Geometric group theory down
under''}\publ  de Gruyter\yr 1999
\endref

\bref\by F. Paulin  \paper Outer automorphisms of hyperbolic groups and small
actions on
$R$-trees,
{\rm pp\. 331--343 in ``Arboreal group theory (R.C. Alperin, ed.)''}\jour
MSRI Publ. 19\publ SpringerVerlag\yr 1991 \endref

\bref  \by M.R. Pettet\paper Virtually free groups with finitely many outer
automorphisms\jour Trans. AMS
 \vol349
\yr1997\pages4565--4587
\endref

\bref \by E. Rips, Z. Sela \paper Structure and rigidity in  hyperbolic
groups I\jour GAFA \vol4\yr1994\pages337--371\endref

\bref \by Z. Sela \paper Structure and rigidity in (Gromov) hyperbolic groups
and discrete groups in rank $1$ Lie groups II\jour GAFA
\vol7\yr1997\pages561--593\endref

\bref\by J. Shor \paper A Scott conjecture for hyperbolic groups \jour
Preprint\yr
1999
\endref

\endRefs

\address   Laboratoire \'Emile Picard, Umr Cnrs 5580,
Universit\'e Paul Sabatier, 31062 Toulouse Cedex 4, France.\endaddress\email
levitt\@picard.ups-tlse.fr{}{}{}{}{}\endemail

\enddocument